\documentclass[12pt]{article}
\usepackage{latexsym}
\usepackage{enumerate}
\usepackage{epsf,epsfig,amsfonts,a4wide}
\usepackage{amsfonts}
\usepackage{url}
\usepackage{amsmath,amssymb,amsthm}
\usepackage{epsf,epsfig,amsfonts,graphicx}
\usepackage{a4wide}

\usepackage[centerlast]{caption}
\usepackage{graphics,amsmath,amssymb,amscd}
\usepackage{amsbsy}
\usepackage{amsthm}
\usepackage{gensymb}
\usepackage{float}
\usepackage{graphicx}

\newcommand{\be}{\begin{equation}}
\newcommand{\ee}{\end{equation}}

\newcommand{\const}{\mathrm{const}}

\usepackage{amsthm}
\usepackage{mathdots}

\usepackage{amscd}

\newcommand{\eps}{\varepsilon}
\newcommand{\ph}{\varphi}
\newcommand{\thet}{\vartheta}
\newcommand{\me}{\mathrm{e}}
\newcommand{\mi}{\mathrm{i}}
\newcommand{\dif}{\mathrm{d}}

\newtheorem{thm}{Theorem}[section]
\newtheorem{lem}{Lemma}[section]

\newtheorem{cor}{Corollary}[section]

\title{On the accuracy of conservation of adiabatic invariants in slow-fast systems}
\author{Tan Su\footnote{ School of Mathematics,
 Loughborough University, UK, LE11 3TU }}

\begin{document}
\maketitle

\begin{abstract}
Let the adiabatic invariant of action variable in slow-fast Hamiltonian system with two degrees of freedom have two limiting values along the trajectories as time tends to infinity. The difference of two limits is exponentially small in analytic systems. An iso-energetic reduction and canonical transformations are applied to transform the slow-fast systems to form of systems depending on slowly varying parameters in a complexified phase space. On the basis of this method an estimate for the accuracy of conservation of adiabatic invariant is given for such systems.
\end{abstract}

\vskip 20pt

\section{Introduction}

Consider a Hamiltonian system with two degrees of freedom. The
Hamiltonian $E$ depends slowly on coordinate $\eps^{-1}x$ and fast
on coordinate $q$:
$$E=E(p,q,y,x)$$
where $q$, $\eps^{-1}x$ are coordinates, and $p$, $y$ are their conjugate momenta.

The variation of the variables $p$, $q$, $y$, $x$ is described by the differential equations

$$
\dot{p}=-\frac{\partial E}{\partial q}, \quad\dot{q}=\frac{\partial E}{\partial p}, \quad
\dot{y}=-\eps\frac{\partial E}{\partial x}, \quad\dot{x}=\eps\frac{\partial E}{\partial y}.
$$

The variables $p$, $q$ are called fast, and $y$, $x$ are called slow variables. This system is called a slow-fast Hamiltonian system. The system with one degree of freedom in which $y=\mathrm{const}$, $x=\mathrm{const}$ is called unperturbed or fast system. Assume that when $y, x=\mathrm{const}$, the phase portrait of the unperturbed system contains a domain filled by closed trajectories \cite{2}:

\begin{figure}[H]
\centering
\includegraphics[width=7cm]{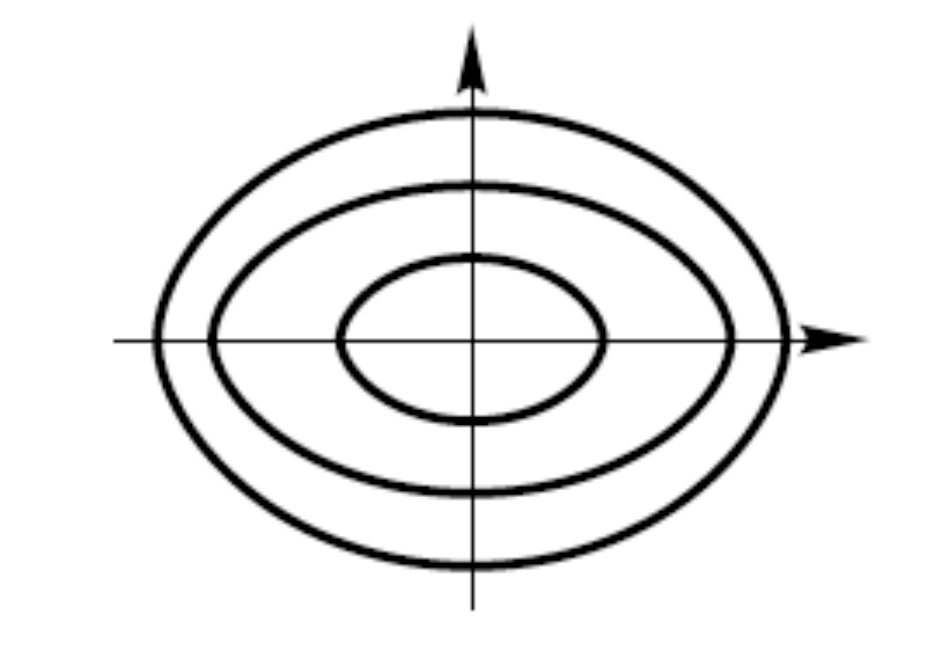}
\caption{}
\end{figure}

Here the frequency of motion is non-zero. So we can introduce action-angle variables

$$I=I(p,q,y,x), \quad \ph=\ph(p,q,y,x)\mod{2\pi}$$

The action $I$ represents the area surrounded by each trajectory. It is equal to this area divided by $2\pi$. Denote $H_0(I,y,x)$ the Hamiltonian $E$ expressed in terms of $I$, $y$, $x$. The approximation in which $I=\mathrm{const}$ and dynamics of $y$, $x$ is described by Hamiltonian system with Hamiltonian $H_0(I,y,x)$ is called an adiabatic approximation. In this approximation,

$$E(p,q,y,x)\equiv H_0(I,y,x),$$
$$\dot{y}=-\eps\frac{\partial H_0(I,y,x)}{\partial x}, \quad\dot{x}=\eps\frac{\partial H_0(I,y,x)}{\partial y}.$$
\vskip 5pt

We will assume that in the phase portrait of this system there is a domain in which along trajectories $|x|\to\infty$, $y\to\const$ as time $t\to\pm\infty$. Thus trajectories have a form shown in one of figures below:

\begin{figure}[H]
\centering
\includegraphics[width=16cm]{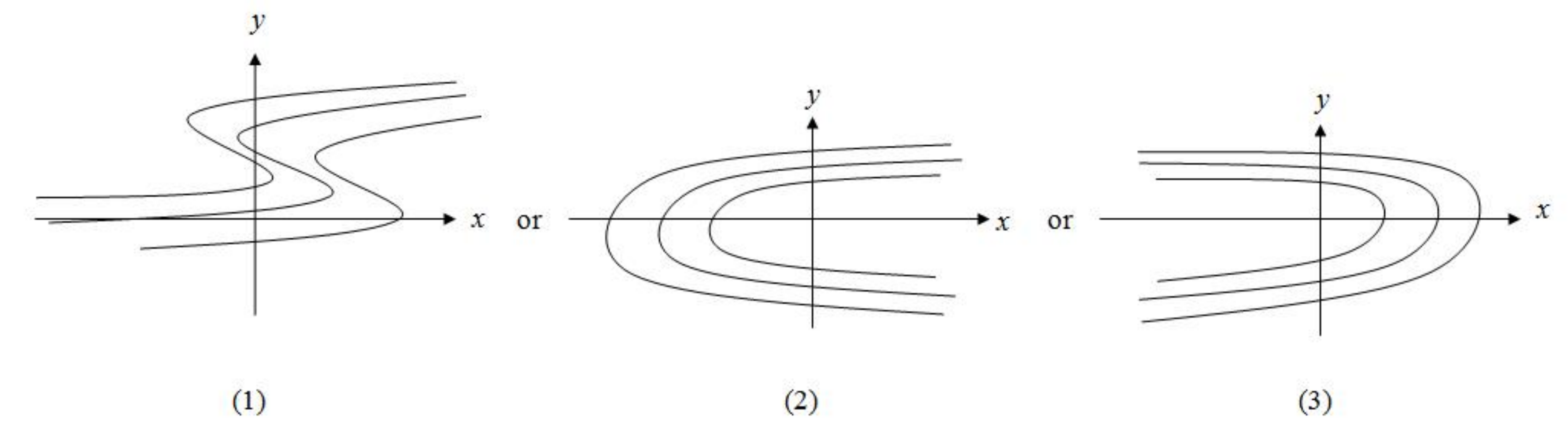}
\caption{}
\end{figure}

We consider the case when the adiabatic invariant of the action $I$ along a trajectory of the original system possesses limiting values $I_{\pm}$ as $t\to\pm\infty$. Their difference $\Delta I=I_+-I_-$ is called an accuracy of persistence of the adiabatic invariant in slow-fast systems.

If this system is analytic, the value of $\Delta I$ is exponentially small \cite{2}:

$$\Delta I=O(\me^{-\frac{\gamma}{\eps}}), \quad \gamma=\mathrm{const}>0.$$

The goal of this paper is to find an estimate of constant $\gamma$. After constructing some assumptions, we will then formulate and prove the conclusion.

\bigskip
\section{Reduction to a standard form}

For the original Hamiltonian $E(p,q,y,x)$ with $y, x=\mathrm{const}$, $q$ is the coordinate and $p$ is conjugate momentum. In the unperturbed system, denote $T$ the period of the trajectory and $\omega$ the frequency. Let $S(I,q,y,x)$ denote the generating function of canonical transformation $(p,q)\mapsto(I,\ph)$. We have the relations \cite{3}:

$$\ph=\frac{\partial S(I,q,y,x)}{\partial I}, \quad p=\frac{\partial S(I,q,y,x)}{\partial q}. $$

For a fixed trajectory $E(p,q,y,x)=h=H_0(I,y,x)$, where $h$ is a constant, we can express $p$ as

$$p=P(h,q,y,x)$$ if $\displaystyle\frac{\partial E}{\partial p}\ne0$. So the generating function has the following form \cite{3}:
$$S(I,q,y,x)=\int\limits_{q_0}^q P(h(I,y,x),q_1,y,x)\,\dif q_1,$$
where the constant $q_0$ is the initial value of $q_1$.

\begin{lem} (see, e.g. \cite{9})
$$\frac{\partial S}{\partial x}=\int\limits_0^t\left(\left<\frac{\partial E}{\partial x}\right>-\frac{\partial E}{\partial x}\right)\,\dif t_1, $$
$$\frac{\partial S}{\partial y}=\int\limits_0^t\left(\left<\frac{\partial E}{\partial y}\right>-\frac{\partial E}{\partial y}\right)\,\dif t_1$$
where $$\left<\frac{\partial E}{\partial x}\right>=\frac 1T\oint\frac{\partial E}{\partial x}\,\dif t_1, \quad \left<\frac{\partial E}{\partial y}\right>=\frac 1T\oint\frac{\partial E}{\partial y}\,\dif t_1.$$
\end{lem}

Now consider another canonical transformation $(p,q,y,x)\mapsto(\bar I, \bar\ph, \bar y, \bar x)$ with generating function $\eps^{-1}\bar y x+S(\bar I, q, \bar y, x)$. The old conjugate variables are $(p,q)$ and $(y, \eps^{-1}x)$, and new pairs after transformation are $(\bar I, \bar\ph)$ and $(\bar y, \eps^{-1}\bar x)$.

\begin{lem} (see, e.g. \cite{2})
$$y=\bar y+O(\eps), \  x=\bar x+O(\eps), \  I=\bar I+O(\eps), \  \ph=\bar\ph+O(\eps)$$
and new Hamiltonian is $$H=H_0(\bar I, \bar y, \bar x)+\eps H_1(\bar I, \bar\ph, \bar y, \bar x,\eps).$$
\end{lem}

We will consider the case that there exist limiting values of action variables $I$ and $\bar I$ as time tends to infinity. The conditions for this and the proof of existence of limits will be added later. Then we will see that the limiting values of $I$ and $\bar I$ are equal.
\vskip 10pt

So far we have got a new Hamiltonian $$H(\bar I, \bar\ph, \bar y, \bar x,\eps)=H_0(\bar I, \bar y, \bar x)+\eps H_1(\bar I, \bar\ph, \bar y, \bar x,\eps).$$
For simplicity, omitting the bar symbols and dependence on $\eps$, we obtain a Hamiltonian system in standard form:
$$H(I, \ph, y,x)=H_0(I, y, x)+\eps H_1(I, \ph, y,x)$$ which has motion
\begin{eqnarray*}
&&\dot I=-\eps\frac{\partial H_1}{\partial\ph}, \quad \dot\ph=\frac{\partial H_0}{\partial I}+\eps\frac{\partial H_1}{\partial I}, \\
&&\dot y=-\eps\left(\frac{\partial H_0}{\partial x}+\eps\frac{\partial H_1}{\partial x}\right), \quad \dot x=\eps\left(\frac{\partial H_0}{\partial y}+\eps\frac{\partial H_1}{\partial y}\right).
\end{eqnarray*}

\bigskip
\section{Statement of the problem}

It is assumed, that in adiabatic approximation, $|x|\to\infty$, $y\to\const$ as $t\to\pm\infty$ (see Figure 2). We will assume also, that when $x$, $y$ are changing in accordance with adiabatic approximation, and $p$, $q$ are fixed:
$$\frac{\partial E(p,q,x,y)}{\partial x}\to0, \quad \text{as\ }t\to\pm\infty.$$

For definiteness we will consider the case when the trajectories in adiabatic approximation have the form shown in Figure 2(1). Result will be valid for the other two cases of Figure 2.

The function $$\omega_0(I,y,x)=\displaystyle\frac{\partial H_0(I,y,x)}{\partial I}$$ is the frequency of unperturbed motion.

Let us assume that the following conditions are fulfilled.

\vskip 2pt
{\bf \underline{Assumption $1^\circ$.}}
\vskip 2pt

The functions $H_0$, $H_1$ can be analytically extended into a complex domain $D=D_I\times D_\ph\times D_{xy}$, where $D_I$ is a neighbourhood of a given real point $I_*$,
\begin{figure}[H]
\centering
\includegraphics[width=8cm]{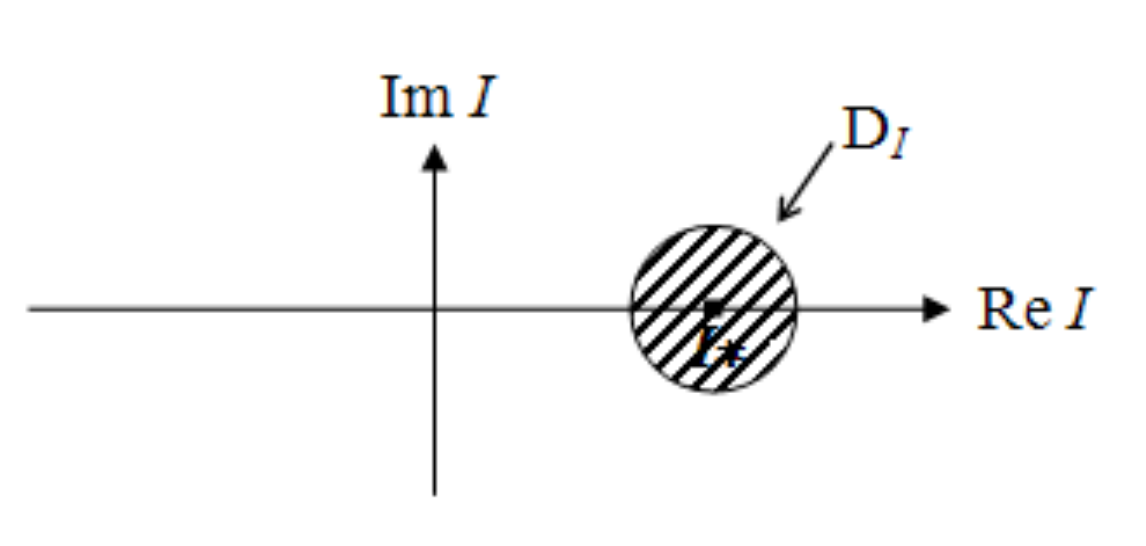}
\caption{}
\end{figure}

\noindent $D_\ph$ is a strip of a fixed width about the real axis,
\begin{figure}[H]
\centering
\includegraphics[width=8cm]{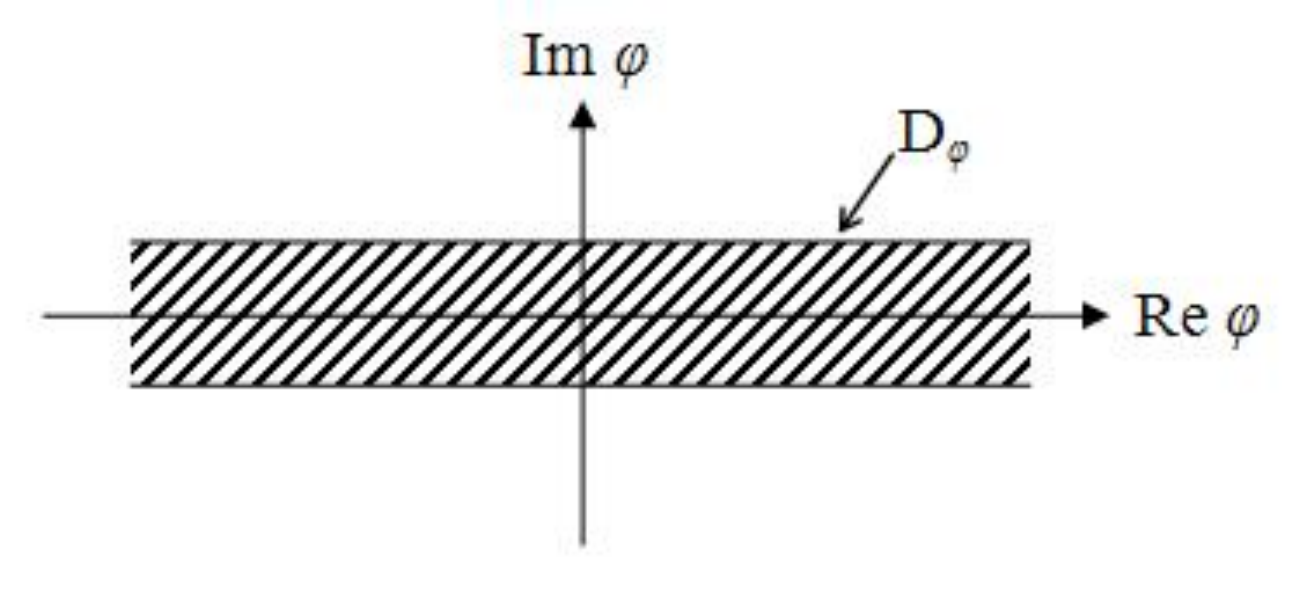}
\caption{}
\end{figure}

\noindent and $D_{xy}$ is some domain in the complex plane $\mathbb C^2$. The function $\omega_0(I,y,x)$ does not vanish in $D$ and $\left|\omega_0\right|>\mathrm{const}$. Function $H_0(I,y,x)$ satisfies
$$\left|\frac{\partial H_0}{\partial I}\right|<\mathrm{const}, \quad \left|\frac{\partial H_0}{\partial y}\right|<\mathrm{const}, \quad \left|\frac{\partial H_0}{\partial x}\right|<\mathrm{const},\quad \left(\frac{\partial H_0}{\partial y}\right)^2+\left(\frac{\partial H_0}{\partial x}\right)^2>\mathrm{const}.$$

Now let us consider approximate Hamiltonian $H_0(I,y,x)$ which has motion
$$\dot y=-\eps\frac{\partial H_0}{\partial x}, \quad\dot x=\eps\frac{\partial H_0}{\partial y}.$$

For a fixed $I_0$, we have $H_0(I_0,y,x)=h=\mathrm{const}$. Assume that curves $H_0=\mathrm{const}$ are not closed:
\begin{figure}[H]
\centering
\includegraphics[width=10cm]{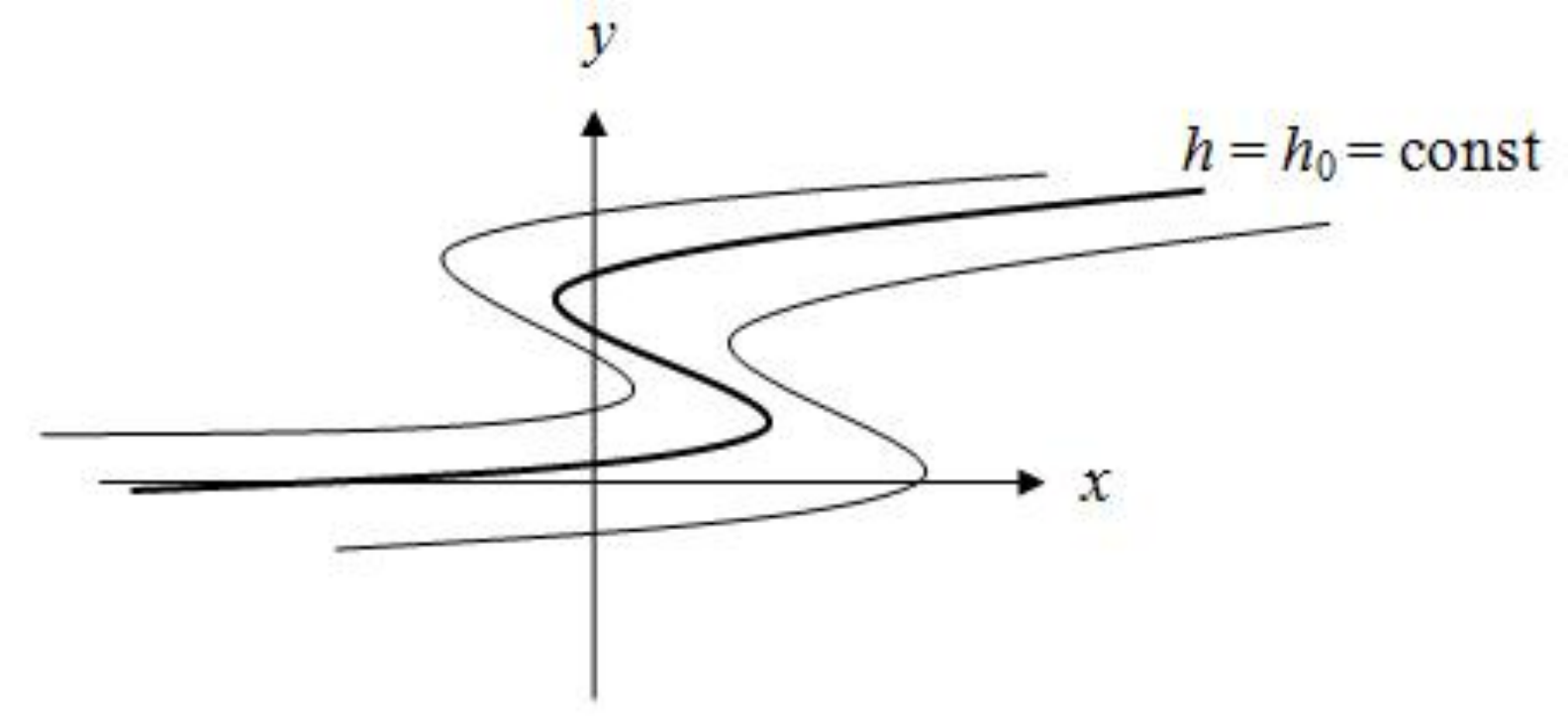}
\caption{}
\end{figure}

\noindent Introduce the slow time $$\tau=\eps t.$$ The differential equations of motion are
$$\frac{\dif y}{\dif\tau}=-\frac{\partial H_0}{\partial x}, \quad\frac{\dif x}{\dif\tau}=\frac{\partial H_0}{\partial y}.$$

Take $h=h_0$, $h_0$ is a point in some interval $D_h$ with centre $h_*$, where $h_*$ is a given real point:
\begin{figure}[H]
\centering
\includegraphics[width=7cm]{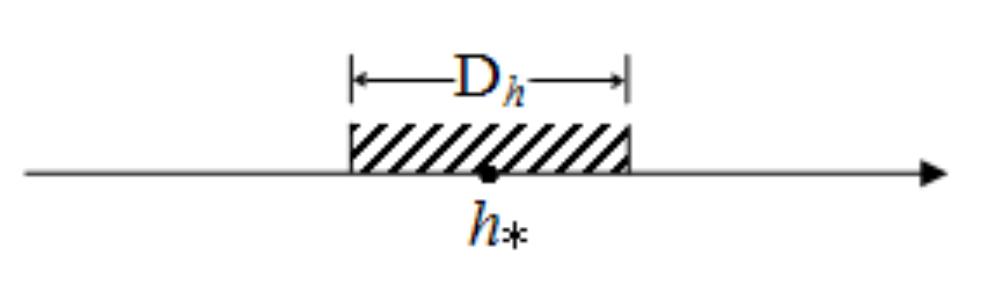}
\caption{}
\end{figure}

\noindent We can find the solution for describing the motion for approximate Hamiltonian $H_0$ with $H_0(I_0,y,x)=h_0$:
$$\left\{ \begin{array}{ccc}
y=Y(\tau, I_0, h_0) \\
x=X(\tau, I_0, h_0)
\end{array}\right.
\qquad I_0\in D_I, \quad h_0\in D_h.$$
From now on we omit the dependence on $h_0$.

\vskip 2pt
{\bf \underline{Assumption $2^\circ$.}}
\vskip 2pt

The solutions $Y(\tau, I_0)$, $X(\tau, I_0)$ can be analytically extended into a strip
$$D_{\tau}=\{\tau\colon|\mathrm{Im}\ \tau|<\sigma+\delta|\mathrm{Re}\ \tau|\}.$$
\begin{figure}[H]
\centering
\includegraphics[width=8cm]{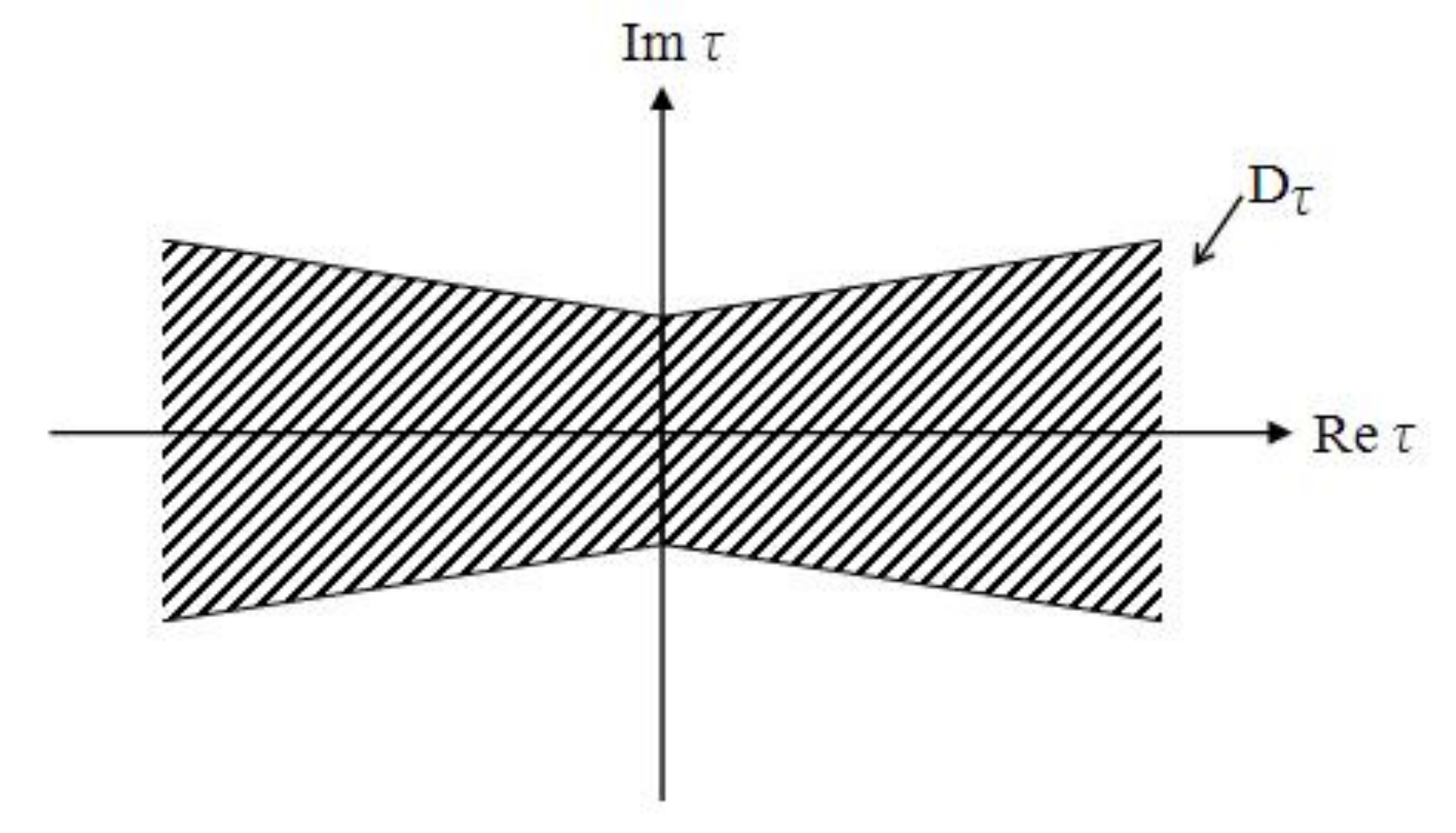}
\caption{}
\end{figure}

We suppose that for any initial point, the solution tends to infinity as time tends to infinity, i.e.
$$\lim\limits_{\mathrm{Re}\,\tau\to\pm\infty}\big|X(\tau,I_0)\big|=\infty.$$

$H_0(I,y,x)$ satisfies
$$\lim\limits_{\mathrm{Re}\,\tau\to\pm\infty}\frac{\partial H_0(I, Y(\tau, I_0), X(\tau, I_0))}{\partial x}=0.$$
$H_1(I,\ph,y,x)$ satisfies
$$|H_1(I, \ph, Y(\tau, I_0), X(\tau, I_0))|<\frac c{1+\left|\int\limits_0^{\tau}\omega_0(I_0, Y(\tau_1, I_0), X(\tau_1, I_0))\,\dif\tau_1\right|^{2+\nu}}.$$
$\omega_0(I,y,x)$ satisfies
$$\mathrm{Im}\ \omega_0(I_0, Y(\tau, I_0), X(\tau, I_0))\rightrightarrows0, \quad \mathrm{as}\ \mathrm{Re}\;\tau\to\pm\infty, \ \mathrm{Im}\;I_0\to0.$$
Here $\sigma$, $\delta$, $c$, $\nu$ are positive constants.

\begin{lem}
For $I\in\widetilde D_I=D_I-\delta_I$, $\ph\in\widetilde D_\ph=D_\ph-\delta_\ph$, and \\
$\tau\in\widetilde D_\tau=\{\tau\colon|\mathrm{Im}\ \tau|<\sigma-\delta_\tau+\delta|\mathrm{Re}\ \tau|\}$, \vskip 1pt
\begin{eqnarray*}
&&\left|\frac{\partial H_1(I, \ph, Y(\tau, I_0), X(\tau, I_0))}{\partial I}\right|<\frac{\mathrm{const}}{1+\left|\int\limits_0^{\tau}\omega_0(I_0, Y(\tau_1, I_0), X(\tau_1, I_0))\,\dif\tau_1\right|^{2+\nu}}, \\
&&\left|\frac{\partial H_1(I, \ph, Y(\tau, I_0), X(\tau, I_0))}{\partial\ph}\right|<\frac{\mathrm{const}}{1+\left|\int\limits_0^{\tau}\omega_0(I_0, Y(\tau_1, I_0), X(\tau_1, I_0))\,\dif\tau_1\right|^{2+\nu}}, \\
&&\left|\frac{\partial H_1(I, \ph, Y(\tau, I_0), X(\tau, I_0))}{\partial\tau}\right|<\frac{\mathrm{const}}{1+\left|\int\limits_0^{\tau}\omega_0(I_0, Y(\tau_1, I_0), X(\tau_1, I_0))\,\dif\tau_1\right|^{2+\nu}},
\end{eqnarray*}
where $\delta_I$, $\delta_\ph$, $\delta_\tau$ are positive constants.
\end{lem}

\begin{proof}
$\mathrm{Im}\ \omega_0(I_0, Y(\tau, I_0), X(\tau, I_0))\rightrightarrows0, \quad \mathrm{as}\ \mathrm{Re}\;\tau\to\pm\infty, \ \mathrm{Im}\;I_0\to0$ \\ implies that
$\forall\,\rho>0,\ \exists\,\Sigma>0,\ \Gamma>0,\ \mathrm{such\ that}\ \mathrm{if}\ |\mathrm{Im}\;I_0|<\Sigma,\ |\mathrm{Re}\;\tau|>\Gamma,\ \mathrm{then}$
$$|\mathrm{Im}\ \omega_0|<\rho.$$
Since we know that $|\omega_0|>c_1=\mathrm{const}$, assume $|\mathrm{Im}\ \omega_0|<\displaystyle\frac{c_1}{10}$ for $|\mathrm{Im}\;I_0|<\Sigma$, $|\mathrm{Re}\;\tau|>\Gamma$, and therefore, $$|\mathrm{Re}\ \omega_0|>\displaystyle\frac9{10}c_1.$$

Let $\tau=\lambda+\mi\mu$. The integration $\left|\int\limits_0^{\tau}\omega_0(I_0, Y(\tau_1, I_0), X(\tau_1, I_0))\,\dif\tau_1\right|$ does not depend on the path.

\vskip 10pt
For $|\mathrm{Re}\;\tau|>\max\Big\{\Gamma,\ 2\Big(\dfrac9{10}c_1-\rho\Big)^{-1}\Big(\dfrac9{10}c_1+\rho+a+b\Big)\Gamma\Big\}$, where $a$, $b$ are positive constants satisfying
$$\int\limits_0^\Gamma|\mathrm{Re}\;\omega_0|\,\dif\lambda\leq a\cdot\Gamma,\quad \int\limits_0^\Gamma|\mathrm{Im}\;\omega_0|\,\dif\lambda\leq b\cdot\Gamma,$$
\vskip-15pt
\begin{figure}[H]
\centering
\includegraphics[width=6cm]{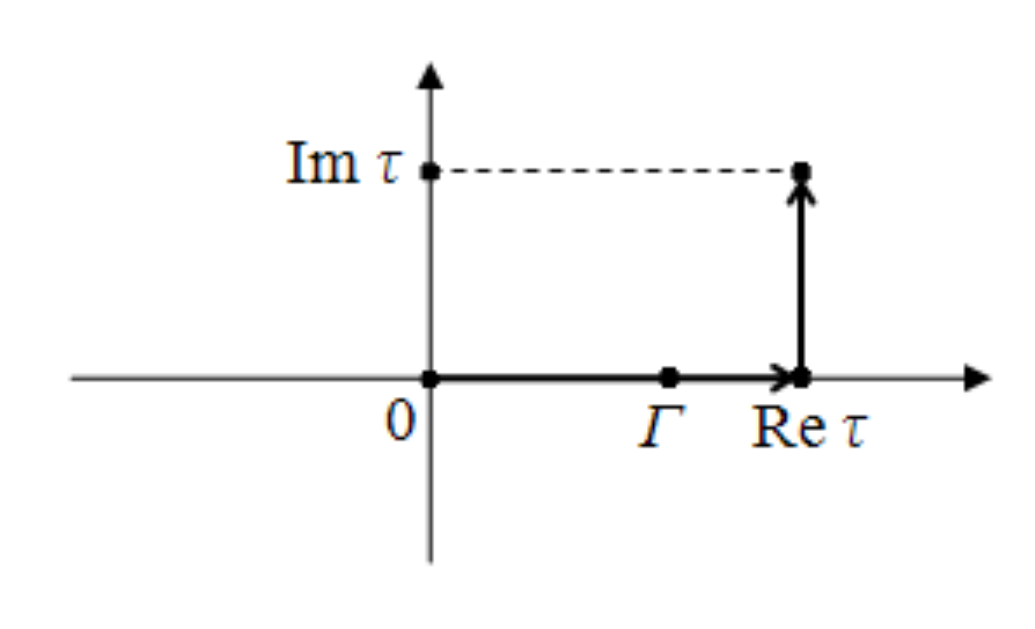}
\caption{}
\end{figure}

we have
\begin{eqnarray*}
\left|\int\limits_0^\tau\omega_0\,\dif\tau_1\right|&=&\left|\int\limits_0^{\mathrm{Re}\,\tau}\omega_0\,\dif\tau_1+\int\limits_{\mathrm{Re}\,\tau}^\tau\omega_0\,\dif\tau_1\right| \\
&=&\left|\int\limits_0^{\mathrm{Re}\,\tau}\omega_0\,\dif (\lambda+\mi\mu)+\int\limits_{\mathrm{Re}\,\tau}^\tau\omega_0\,\dif (\lambda+\mi\mu)\right| \\
&=&\left|\int\limits_0^{\mathrm{Re}\,\tau}\omega_0\,\dif\lambda+\int\limits_{\mathrm{Re}\,\tau}^\tau\omega_0\,\dif(\mi\mu)\right| \\
&=&\left|\int\limits_0^{\mathrm{Re}\,\tau}\mathrm{Re}\;\omega_0\,\dif\lambda+\int\limits_0^{\mathrm{Re}\,\tau}\mi\;\mathrm{Im}\;\omega_0\,\dif\lambda+\int\limits_{\mathrm{Re}\,\tau}^\tau\mathrm{Re}\;\omega_0\,\dif(\mi\mu)+\int\limits_{\mathrm{Re}\,\tau}^\tau\mi\;\mathrm{Im}\;\omega_0\,\dif(\mi\mu)\right| \\
&=&\left|\int\limits_0^\Gamma\mathrm{Re}\;\omega_0\,\dif\lambda+\int\limits_\Gamma^{\mathrm{Re}\,\tau}\mathrm{Re}\;\omega_0\,\dif\lambda+\int\limits_0^\Gamma\mi\;\mathrm{Im}\;\omega_0\,\dif\lambda+\int\limits_\Gamma^{\mathrm{Re}\,\tau}\mi\;\mathrm{Im}\;\omega_0\,\dif\lambda \right. \\
&&\left. +\int\limits_{\mathrm{Re}\,\tau}^\tau\mathrm{Re}\;\omega_0\,\dif(\mi\mu)+\int\limits_{\mathrm{Re}\,\tau}^\tau\mi\;\mathrm{Im}\;\omega_0\,\dif(\mi\mu)\right| \\
\end{eqnarray*}
\begin{eqnarray*}
&\geq&\left|\left|\int\limits_\Gamma^{\mathrm{Re}\,\tau}\mathrm{Re}\;\omega_0\,\dif\lambda+\int\limits_{\mathrm{Re}\,\tau}^\tau\mathrm{Re}\;\omega_0\,\dif(\mi\mu)\right|-\left|\int\limits_\Gamma^{\mathrm{Re}\,\tau}|\mathrm{Im}\;\omega_0|\,\dif\lambda+\int\limits_{\mathrm{Re}\,\tau}^\tau|\mathrm{Im}\;\omega_0|\,\dif(\mi\mu) \right.\right.\\
&&\left.\left. +\int\limits_0^\Gamma|\mathrm{Re}\;\omega_0|\,\dif\lambda+\int\limits_0^\Gamma|\mathrm{Im}\;\omega_0|\,\dif\lambda\right|\right| \\
&>&\left|\Big|\frac9{10}c_1\cdot(\mathrm{Re}\;\tau-\Gamma)+\frac9{10}c_1\cdot\mi\;\mathrm{Im}\;\tau\Big|-\Big|\rho\cdot(\mathrm{Re}\;\tau-\Gamma)+\rho\cdot\mi\;\mathrm{Im}\;\tau+a\Gamma+b\Gamma\Big|\right| \\
&=&\left|\Big|\frac9{10}c_1\cdot\tau-\frac9{10}c_1\cdot\Gamma\Big|-\Big|\rho\cdot\tau-\rho\cdot\Gamma+a\Gamma+b\Gamma\Big|\right| \\
&\geq&\left|\Big(\frac9{10}c_1|\tau|-\frac9{10}c_1\Gamma\Big)-\Big(\rho|\tau|+\rho\Gamma+a\Gamma+b\Gamma\Big)\right| \\
&=&\left|\Big(\frac9{10}c_1-\rho\Big)|\tau|-\Big(\frac9{10}c_1+\rho+a+b\Big)\Gamma\right| \\
&>&\Big(\frac9{10}c_1-\rho\Big)|\tau|
\end{eqnarray*}
as $\Big(\dfrac9{10}c_1-\rho\Big)|\tau|>2\Big(\dfrac9{10}c_1+\rho+a+b\Big)\Gamma$.

We can find a positive constant $c_2<\displaystyle\frac9{10}c_1-\rho$, such that
$$\left|\int\limits_0^{\tau}\omega_0(I_0, Y(\tau_1, I_0), X(\tau_1, I_0))\,\dif\tau_1\right|>\left(\frac9{10}c_1-\rho\right)|\tau|>c_2|\tau|.$$

For $|\mathrm{Re}\;\tau|\leq\max\Big\{\Gamma,\ 2\Big(\dfrac9{10}c_1-\rho\Big)^{-1}\Big(\dfrac9{10}c_1+\rho+a+b\Big)\Gamma\Big\}$, it is evident that we can find another positive constant $c_3$ such that
$$\left|\int\limits_0^{\tau}\omega_0(I_0, Y(\tau_1, I_0), X(\tau_1, I_0))\,\dif\tau_1\right|>c_3|\tau|.$$
So
\begin{eqnarray*}
|H_1(I, \ph, Y(\tau, I_0), X(\tau, I_0))|&<&\frac c{1+\left|\int\limits_0^{\tau}\omega_0(I_0, Y(\tau_1, I_0), X(\tau_1, I_0))\,\dif\tau_1\right|^{2+\nu}} \\
&<&\frac c{1+(\const\cdot|\tau|)^{2+\nu}} \\
&<&\frac{\const}{1+|\tau|^{2+\nu}}.
\end{eqnarray*}

From Cauchy estimate \cite{8}, we can easily get that, in $\widetilde D_I$ and $\widetilde D_\ph$,
$$\left|\frac{\partial H_1(I, \ph, Y(\tau, I_0), X(\tau, I_0))}{\partial I}\right|<\frac{\mathrm{const}}{1+|\tau|^{2+\nu}},$$
$$\left|\frac{\partial H_1(I, \ph, Y(\tau, I_0), X(\tau, I_0))}{\partial\ph}\right|<\frac{\mathrm{const}}{1+|\tau|^{2+\nu}}.$$

Now let us prove $\left|\displaystyle\frac{\partial H_1(I, \ph, Y(\tau, I_0), X(\tau, I_0))}{\partial\tau}\right|<\displaystyle\frac{\mathrm{const}}{1+|\tau|^{2+\nu}}$. From
$$D_\tau=\{\tau\colon|\mathrm{Im}\ \tau|<\sigma+\delta|\mathrm{Re}\ \tau|\},$$
for simplicity, assume $\delta_\tau=\displaystyle\frac\sigma{10}$,  and then
$$\widetilde D_\tau=\{\tau\colon|\mathrm{Im}\ \tau|<\sigma-\frac\sigma{10}+\delta|\mathrm{Re}\ \tau|\}.$$
In the following figure, mark the boundary of $D_\tau$ with solid lines, and that of $\widetilde D_\tau$ with dashed lines:
\begin{figure}[H]
\centering
\includegraphics[width=10cm]{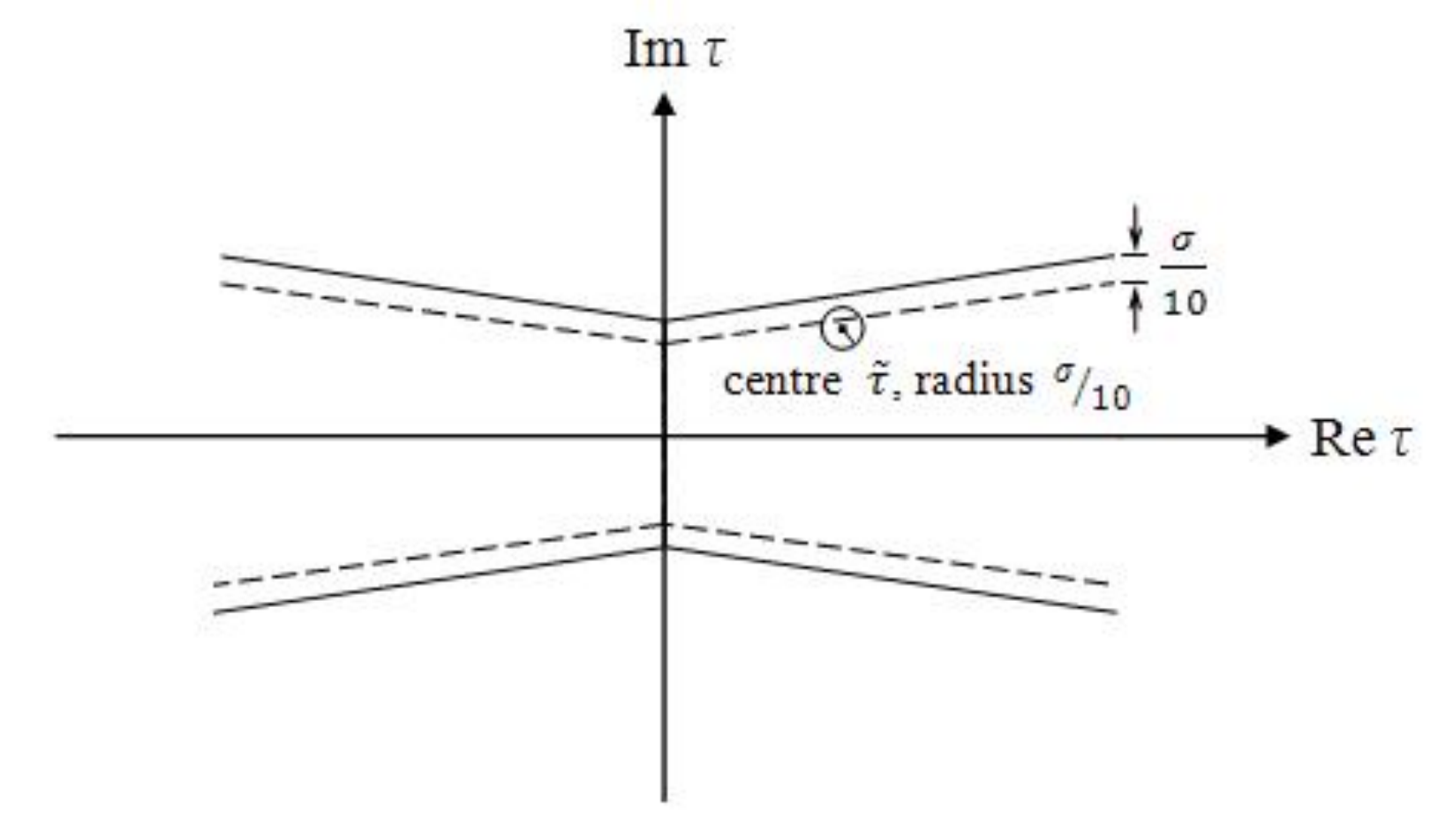}
\caption{}
\end{figure}

For $\widetilde\tau\in\widetilde D_\tau$, the estimate of $|H_1|$ on the circle of radius $\displaystyle\frac\sigma{10}$ around $\widetilde\tau$ is \\ $|H_1|<\displaystyle\frac{\mathrm{const}}{1+\min\limits_\theta\left|\widetilde\tau+\frac\sigma{10}\me^{\mi\theta}\right|^{2+\nu}}$.
If $|\widetilde\tau|>\sigma$, then
$$|H_1|<\frac{\mathrm{const}}{1+\min\limits_\theta\left|\widetilde\tau+\displaystyle\frac\sigma{10}\me^{\mi\theta}\right|^{2+\nu}}<\frac{\mathrm{const}}{1+\left|\displaystyle\frac9{10}\widetilde\tau\right|^{2+\nu}}.$$
So, because of the Cauchy estimate $$\left.\left|\frac{\partial H_1}{\partial\tau}\right|\right|_{\widetilde\tau}<\left.\frac{\mathrm{const}}{1+\left|\frac9{10}\widetilde\tau\right|^{2+\nu}}\right/\!\!\frac\sigma{10}<\frac{\mathrm{const}}{1+|\widetilde\tau|^{2+\nu}}.$$
If $|\widetilde\tau|\leq\sigma$, then because of the Cauchy estimate,
$$\left.\left|\frac{\partial H_1}{\partial\tau}\right|\right|_{\widetilde\tau}<\mathrm{const}\left/\frac\sigma{10}\right.<\frac{\mathrm{const}}{1+|\widetilde\tau|^{2+\nu}}.$$
Therefore, $$\left|\frac{\partial H_1(I, \ph, Y(\tau, I_0), X(\tau, I_0))}{\partial\tau}\right|<\frac{\mathrm{const}}{1+|\tau|^{2+\nu}}.$$

\vskip 10pt
We can find a constant $c_4$ such that $$\left|\int\limits_0^{\tau}\omega_0(I_0, Y(\tau_1, I_0), X(\tau_1, I_0))\,\dif\tau_1\right|<c_4|\tau|.$$ Then
\begin{eqnarray*}
\frac{\mathrm{const}}{1+|\tau|^{2+\nu}}&<&\frac{\mathrm{const}}{1+\left(c_4^{-1}\left|\int\limits_0^{\tau}\omega_0(I_0, Y(\tau_1, I_0), X(\tau_1, I_0))\,\dif\tau_1\right|\right)^{2+\nu}} \\
&<&\frac{\mathrm{const}}{1+\left|\int\limits_0^{\tau}\omega_0(I_0, Y(\tau_1, I_0), X(\tau_1, I_0))\,\dif\tau_1\right|^{2+\nu}}.
\end{eqnarray*}

Therefore, for $I\in\widetilde D_I$, $\ph\in\widetilde D_\ph$, $\tau\in\widetilde D_\tau$,
\begin{eqnarray*}
&&\left|\frac{\partial H_1(I, \ph, Y(\tau, I_0), X(\tau, I_0))}{\partial I}\right|<\frac{\mathrm{const}}{1+\left|\int\limits_0^{\tau}\omega_0(I_0, Y(\tau_1, I_0), X(\tau_1, I_0))\,\dif\tau_1\right|^{2+\nu}}, \\
&&\left|\frac{\partial H_1(I, \ph, Y(\tau, I_0), X(\tau, I_0))}{\partial\ph}\right|<\frac{\mathrm{const}}{1+\left|\int\limits_0^{\tau}\omega_0(I_0, Y(\tau_1, I_0), X(\tau_1, I_0))\,\dif\tau_1\right|^{2+\nu}}, \\
&&\left|\frac{\partial H_1(I, \ph, Y(\tau, I_0), X(\tau, I_0))}{\partial\tau}\right|<\frac{\mathrm{const}}{1+\left|\int\limits_0^{\tau}\omega_0(I_0, Y(\tau_1, I_0), X(\tau_1, I_0))\,\dif\tau_1\right|^{2+\nu}}.
\end{eqnarray*}

\end{proof}

\vskip 2pt
{\bf \underline{Assumption $3^\circ$.}}
\vskip 2pt

The level lines
$$\mathrm{Im}\int\limits_0^\tau\omega_0(I,Y(\tau_1,I),X(\tau_1,I))\,\dif\tau_1=B=\mathrm{const}, \quad 0\leq|B|\leq\gamma, \quad I\in D_I$$
lie in the domain $D_\tau$ and have a positive distance from the boundary of $D_\tau$.

\vskip 10pt
Now consider the exact solution $I(t)$, $\ph(t)$, $y(t)$, $x(t)$ of the Hamiltonian system with Hamiltonian
$$H(I, \ph, y,x)=H_0(I, y, x)+\eps H_1(I, \ph, y,x)$$ with real initial conditions $I(0)$, $\ph(0)$, $y(0)$, $x(0)$ in $D=D_I\times D_\ph\times D_{xy}$ at \\ $t=\tau/\eps=0$, and $H(I(0),\ph(0),y(0),x(0))=h_0$.
The adiabatic invariant of the action $I(t)$ satisfies
$$\frac{\dif I}{\dif t}=-\eps\frac{\partial H_1(I,\ph,y,x)}{\partial\ph}$$
and thus
$$I(t)=I(0)-\eps\int\limits_0^t\frac{\partial H_1(I(t_1),\ph(t_1),y(t_1),x(t_1))}{\partial\ph}\,\dif t_1.$$

\begin{thm}
There exist limiting values $$I_\pm=\lim\limits_{t\to\pm\infty}I(t)$$
and their difference $\Delta I=I_+-I_-$ satisfies the estimate
$$\Delta I=O(\me^{-\frac\gamma\eps}), \quad \gamma=\mathrm{const}>0$$
with the constant $\gamma$ introduced in Assumption \textup{3}$^\circ$.
\end{thm}

{\bf Remark.} {The method of continuous averaging \cite{10} gives the same estimate for $\gamma$, as the above-stated theorem.}

\section{Example}

Paper \cite{5} gives us several examples of slow-fast Hamiltonian systems, in one of which the Hamiltonian is as follows:
$$H(I,\ph,y,x)=\omega I+\frac{y^2}2+V_0(x)+\eps g(\ph)V_1(x).$$
Here $y$, $x$ are slow variables and $I$, $\ph$ are fast ones. We explicitly treat the case of $V_0(x)=\me^{-x}$. Concerning $V_1$, we need it to be real analytic and tend to zero as $x\to+\infty$. For simplicity we choose $V_1\equiv\me^{-x}$. The unperturbed Hamiltonian system is
$$H_0(I,y,x)=\omega I+\frac{y^2}2+V_0(x)=\omega I+\frac{y^2}2+\me^{-x}.$$

Let $H_0(I,y,x)=\omega I+\dfrac{y^2}2+\me^{-x}=h_0=\const$ be the total energy of the unperturbed system. We can draw the phase portrait of the system to describe the motion:

\begin{figure}[H]
\centering
\includegraphics[width=9cm]{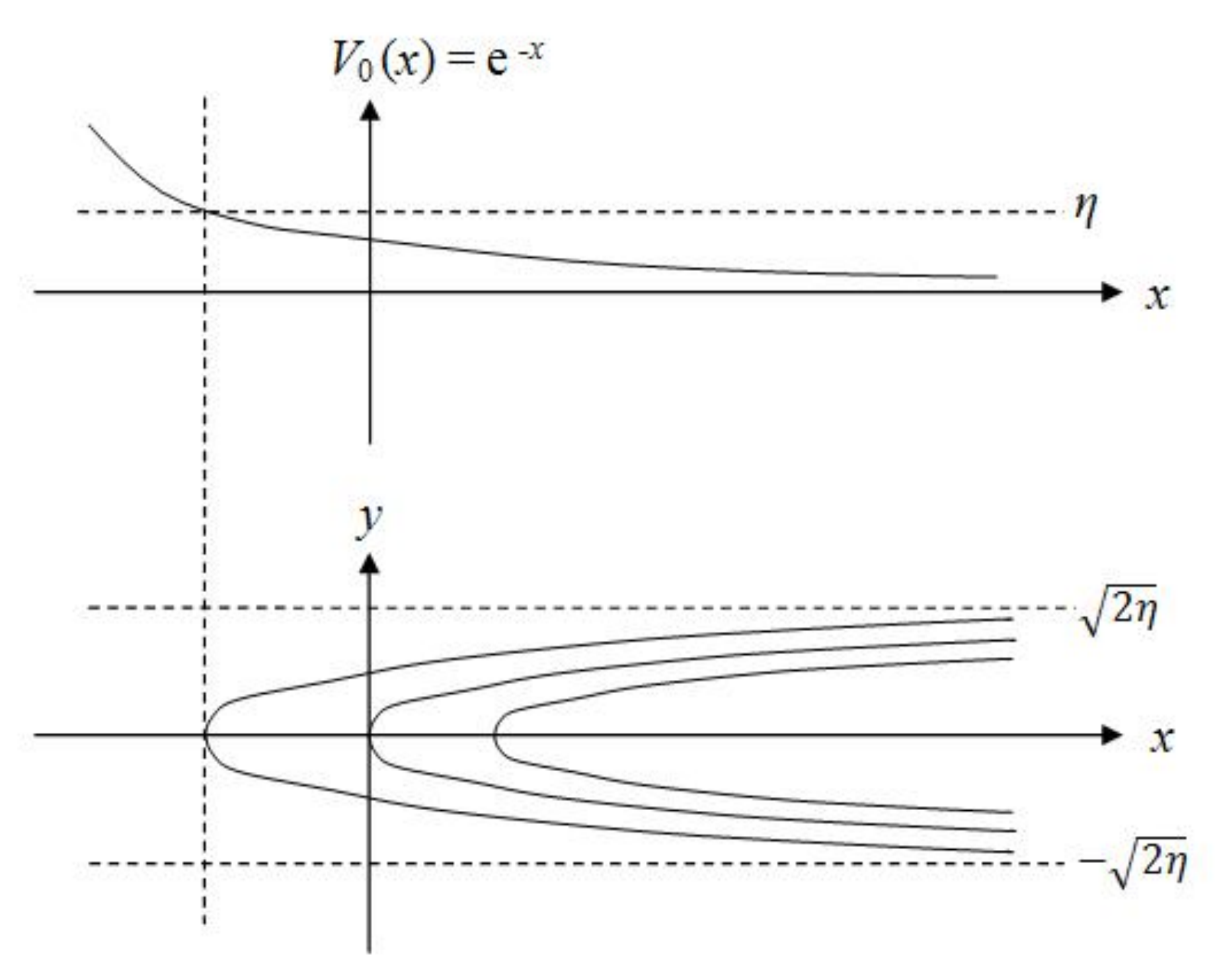}
\caption{}
\end{figure}

Take $\tilde\xi$ as the slow time of motion. From the Hamiltonian $H_0(I,y,x)=\omega I+\dfrac{y^2}2+\me^{-x}=h_0$, we know that
$$\frac{\dif x}{\dif\tilde\xi}=\frac{\partial H_0}{\partial y}=y.$$
So $$\frac12\left(\frac{\dif x}{\dif\tilde\xi}\right)^2+\me^{-x}=h_0-\omega I=\tilde\eta.$$
Thus we can express $x$ and $y$ via $(\tilde\xi, \tilde\eta)$:
$$x=x(\tilde\xi,\tilde\eta),\quad y=y(\tilde\xi,\tilde\eta)$$
with initial data $x^0=x(0,\tilde\eta)$ and $y^0=\sqrt{2(\tilde\eta-\me^{-x^0})}$.

After solving differential equations, we can obtain the solutions as (see \cite{5})
$$x(\tilde\xi,\tilde\eta)=\log\left(\frac1{\tilde\eta}\Big(\cosh\sqrt{\tilde\eta/2}\,\big(\tilde\xi-\tilde\xi^0\big)\Big)^2\right), \quad y(\tilde\xi,\tilde\eta)=\sqrt{2\tilde\eta}\,\tanh\sqrt{\tilde\eta/2}\,\big(\tilde\xi-\tilde\xi^0\big).$$
For simplicity we take the initial value of $\tilde\xi$ as $\tilde\xi^0=0$. So $x^0=x(0,\tilde\eta)=\log\dfrac1{\tilde\eta}$, $y^0=0$.

Now after the canonical transformation of variables from $(x,y)$ to $(\tilde\xi,\tilde\eta)$, the Hamiltonian $H_0$ becomes $K_0(\tilde\xi,\tilde\eta)=\omega I+\tilde\eta$ and the new Hamiltonian is
$$K(\tilde\xi,\tilde\eta,I,\ph)=\omega I+\tilde\eta+\eps f(\tilde\xi,\tilde\eta)g(\ph)$$
with $$f(\tilde\xi,\tilde\eta)=V_1(x(\tilde\xi,\tilde\eta))=\me^{-\log\left(\frac1{\tilde\eta}\left(\cosh\sqrt{\tilde\eta/2}\,\tilde\xi\right)^2\right)}=\frac{\tilde\eta}{\left(\cosh\sqrt{\tilde\eta/2}\,\tilde\xi\right)^2}.$$
The differential equations for describing the motion are
\begin{eqnarray*}
&&\dot{\tilde\xi}=\eps+\eps^2\frac{\partial f(\tilde\xi,\tilde\eta)}{\partial\tilde\eta}g(\ph), \\
&&\dot{\tilde\eta}=-\eps^2\frac{\partial f(\tilde\xi,\tilde\eta)}{\partial\tilde\xi}g(\ph), \\
&&\dot I=-\eps f(\tilde\xi,\tilde\eta)\frac{\dif g(\ph)}{\dif\ph}, \\
&&\dot\ph=\omega.
\end{eqnarray*}

For $g(\ph)$, we suppose that $g(\ph)$ is analytic in the strip $|\text{Im}\,\ph|<\rho$ with constant $\rho>0$, and bounded:
$$|g(\ph)|\le1, \quad \text{for}\ |\text{Im}\,\ph|<\rho.$$

For $f(\tilde\xi,\tilde\eta)$, we can easily prove that
$$|f(\tilde\xi,\tilde\eta)|<\frac{\const}{1+|\tilde\xi|^{2+\nu}},$$
and $$\left|\frac{\partial f(\tilde\xi,\tilde\eta)}{\partial \tilde\xi}\right|<\frac{\const}{1+|\tilde\xi|^{2+\nu}},$$
$$\left|\frac{\partial f(\tilde\xi,\tilde\eta)}{\partial \tilde\eta}\right|<\frac{\const}{1+|\tilde\xi|^{2+\nu}}.$$

Take $D_I$ as a neighbourhood of a given real point $I_*$, and $D_\ph$ is a strip $|\text{Im}\,\ph|<\rho$. Now let us determine the domain of $\tilde\xi$ by finding singularities of $f(\tilde\xi,\tilde\eta)$:
$$f(\tilde\xi,\tilde\eta)=\frac{\tilde\eta}{\left(\cosh\sqrt{\tilde\eta/2}\,\tilde\xi\right)^2} =\frac{4\tilde\eta}{\left(\me^{\sqrt{\tilde\eta/2}\,\tilde\xi}+\me^{-\sqrt{\tilde\eta/2}\,\tilde\xi}\right)^2}.$$
Points of singularities should satisfy
\begin{eqnarray*}
&&\me^{\sqrt{\tilde\eta/2}\,\tilde\xi}+\me^{-\sqrt{\tilde\eta/2}\,\tilde\xi}=0 \\
&&\me^{2\sqrt{\tilde\eta/2}\,\tilde\xi}+1=0 \\
&&\me^{\sqrt{2\tilde\eta}\,\tilde\xi}=-1 \\
&&\me^{\text{Re}\,\sqrt{2\tilde\eta}\,\tilde\xi}\,\left(\cos\text{Im}\,\sqrt{2\tilde\eta}\,\tilde\xi+\mi\sin\text{Im}\,\sqrt{2\tilde\eta}\,\tilde\xi\right)=-1
\end{eqnarray*}
\begin{eqnarray*}
\text{So}  &&\cos\text{Im}\,\sqrt{2\tilde\eta}\,\tilde\xi=-\me^{-\text{Re}\,\sqrt{2\tilde\eta}\,\tilde\xi}<0, \\
&&\sin\text{Im}\,\sqrt{2\tilde\eta}\,\tilde\xi=0.\\
\text{Therefore,} &&\text{Im}\,\sqrt{2\tilde\eta}\,\tilde\xi=(2k+1)\pi,\quad k\in\mathbb Z,\\
\text{and} &&\text{Re}\,\sqrt{2\tilde\eta}\,\tilde\xi=0.
\end{eqnarray*}
Thus $$\sqrt{2\tilde\eta}\,\tilde\xi=(2k+1)\pi\mi,\quad k\in\mathbb Z.$$
The domain of $\tilde\eta$ denoted as $D_{\tilde\eta}$ is a complex ball around a real value $\tilde\eta_*$. Any value $\tilde\eta^0\in D_{\tilde\eta}$ can be taken to determine the singularities. Therefore, we obtain the singularities as
$$\tilde\xi=(2k+1)\frac\pi{\sqrt{2\tilde\eta^0}\;}\,\mi,\quad k\in\mathbb Z.$$
From the above conclusion, if we take $\tilde\xi$ in the strip
$$D_{\tilde\xi}=\{\tilde\xi\colon|\text{Im}\,\tilde\xi|<\frac\pi{\sqrt{2\tilde\eta^0}\;}+\delta|\text{Re}\,\tilde\xi|\},$$
where $\delta$ is a positive constant, $f(\tilde\xi,\tilde\eta)$ is analytic for any $\tilde\eta\in D_{\tilde\eta}$.

The following figure shows the boundary of domain $D_{\tilde\xi}$:
\begin{figure}[H]
\centering
\includegraphics[width=11cm]{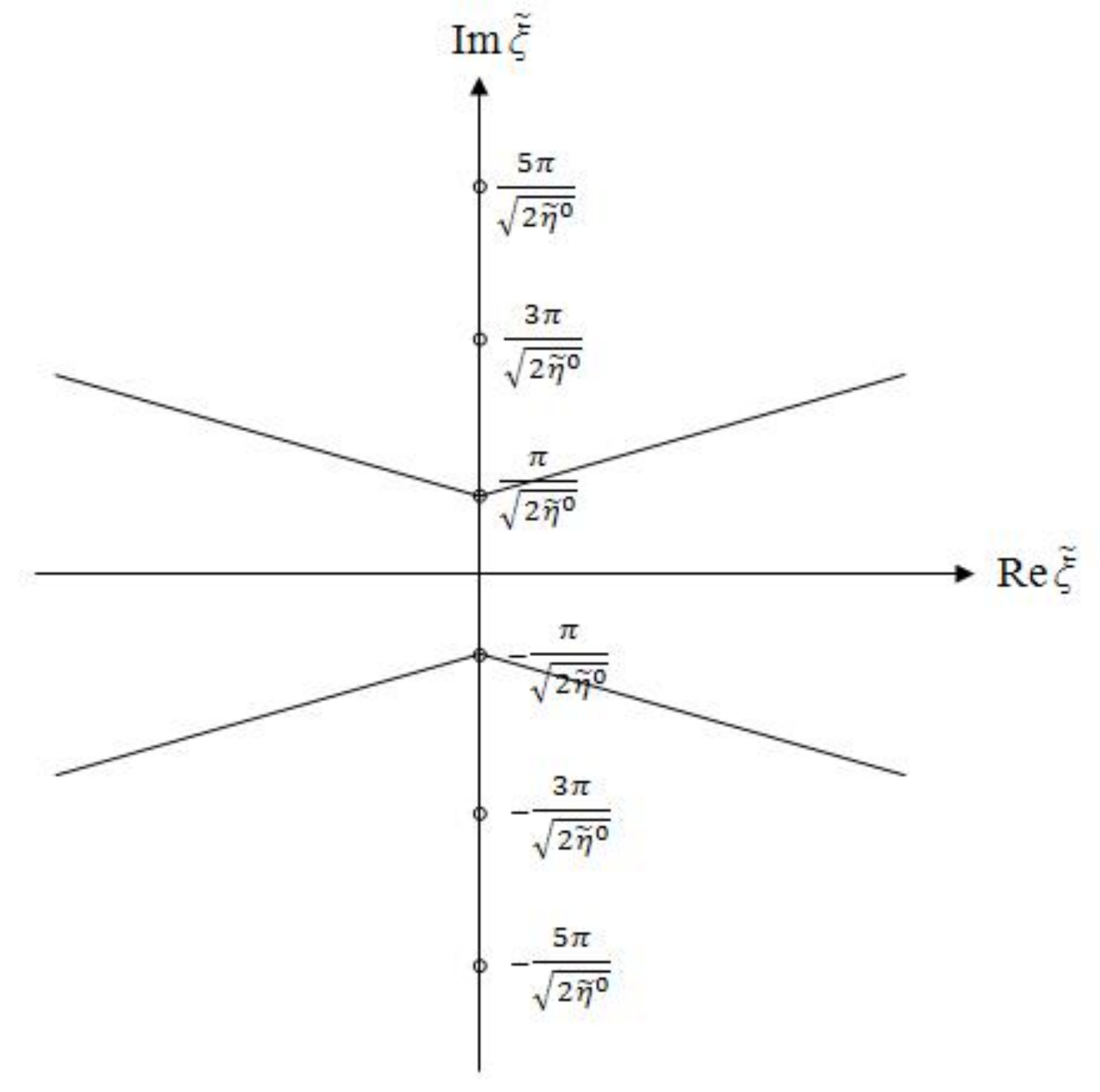}
\caption{}
\end{figure}

We can find the level lines
$$\text{Im}\,\int\limits_0^{\tilde\xi}\omega\,\dif\tilde\xi_1=\omega\text{Im}\,\tilde\xi=B=\const,\quad 0\le|B|\le\gamma$$
lying in the domain $D_{\tilde\xi}$, which have a positive distance from the boundary of this domain. So $\gamma$ has been found as any constant such that
$$\gamma<\frac{\pi\omega}{\sqrt{2\tilde\eta^0}},\quad \rm{for\ such\ definition\ of\ \tilde\eta_0}.$$
According to Theorem 1, the accuracy of conservation of adiabatic invariant $\Delta I$ has the estimate:
$$\Delta I=O(\me^{-\frac\gamma\eps}).$$

This estimate coincides with that in \cite{5}, where it was obtained by a different method.

\section{Proof of the Theorem}

The original system is described by Hamiltonian
$$H(I, \ph, y,x)=H_0(I, y, x)+\eps H_1(I, \ph, y,x)$$ with motion in time $t$:
\begin{eqnarray*}
&&\dot I=-\eps\frac{\partial H_1}{\partial\ph}, \quad \dot\ph=\frac{\partial H_0}{\partial I}+\eps\frac{\partial H_1}{\partial I}, \\
&&\dot y=-\eps\left(\frac{\partial H_0}{\partial x}+\eps\frac{\partial H_1}{\partial x}\right), \quad \dot x=\eps\left(\frac{\partial H_0}{\partial y}+\eps\frac{\partial H_1}{\partial y}\right).
\end{eqnarray*}

As the frequency $\omega_0(I,y,x)=\displaystyle\frac{\partial H_0(I,y,x)}{\partial I}$ does not vanish in our domain $D$, we can express $I$ from $H$ via $\ph$, $y$, $x$ and a constant of energy. This is the procedure of isoenergetic reduction on the energy level $H(I,\ph,y,x)=h_0,\ h_0=\mathrm{const}$:
$$-I(y,x,\ph,h_0)=F_0(y,x,h_0)+\eps F_1(y,x,\ph,h_0).$$
Here $(-I)$ is the new Hamiltonian, $\ph$ is the new time, and $(y, x)$ is in a small neighbourhood of $\widetilde D_{xy}$, which can be written as $\widetilde D_{xy}+\delta_{xy}$, where
$$\widetilde D_{xy}=\bigg\{(y,x)\colon \Big\{ \begin{array}{c}
y=Y(\tau, I_0, h_0)\\
x=X(\tau, I_0, h_0)
\end{array}\!,\ \tau\in D_\tau,\ I_0\in D_I\bigg\},$$
and $\delta_{xy}$ is some positive constant.

We do not indicate dependence of $Y$, $X$ on $\eps$. The differential equations of motion are
$$\frac{\dif y}{\dif\ph}=-\eps\left(\frac{\partial F_0}{\partial x}+\eps\frac{\partial F_1}{\partial x}\right),\quad \frac{\dif x}{\dif\ph}=\eps\left(\frac{\partial F_0}{\partial y}+\eps\frac{\partial F_1}{\partial y}\right).$$
Since we know functions $H_0$ and $H_1$ are analytic in $D$, from the implicit function theorem, functions $F_0$ and $F_1$ are also analytic in $D_\ph\times(\widetilde D_{xy}+\delta_{xy})$.

Now consider the approximate system with Hamiltonian $$F_0(y,x,h_0)=-I_0.$$
Introduce slow time $\xi=\eps\ph$, and the motion is
$$\frac{\dif y}{\dif\xi}=-\frac{\partial F_0}{\partial x},\quad \frac{\dif x}{\dif\xi}=\frac{\partial F_0}{\partial y}.$$
We can find the solution for the above approximate system:
$$\left\{ \begin{array}{ccc}
y=\widehat Y(\xi, I_0, h_0) \\
x=\widehat X(\xi, I_0, h_0)
\end{array}\right.
\qquad I_0\in D_I, \quad h_0\in D_h.$$
In what follows, we do not indicate the dependence on  constant $h_0$. The relation between solutions $Y(\tau,I_0)$, $X(\tau,I_0)$ and $\widehat Y(\xi,I_0)$, $\widehat X(\xi,I_0)$ follows from the formula
$$\xi=\int\limits_0^{\tau}\omega_0(I_0, Y(\tau_1, I_0), X(\tau_1, I_0))\,\dif\tau_1.$$
Solution $\widehat Y(\xi,I_0)$, $\widehat X(\xi,I_0)$ can be analytically extended into some domain $D_\xi$ with respect to $\xi$. The domain $D_\xi$ will be introduced in later lemma. As we know,
$$\mathrm{Im}\ \omega_0(I_0, Y(\tau, I_0), X(\tau, I_0))\rightrightarrows0, \quad \mathrm{as}\ \mathrm{Re}\;\tau\to\pm\infty, \ \mathrm{Im}\;I_0\to0.$$
It implies that $\omega_0$ uniformly tends to a real function depending on Re$\,I_0$ as time tends to infinity and $I_0$ tends to real. Let us denote the limits as $\omega_\pm(\mathrm{Re}\,I_0)$:
$$\omega_0(I_0, Y(\tau_1, I_0), X(\tau_1, I_0))\rightrightarrows\omega_\pm(\mathrm{Re}\,I_0), \quad \mathrm{as} \ \mathrm{Re}\,\tau\to\pm\infty,\ \mathrm{Im}\,I_0\to0.$$

\begin{lem}
For $I_0\in D_I$ close to the imaginary axis, the image of the domain $D_\tau$ under the map $\tau\mapsto\xi$ given by the formula $$\xi=\int\limits_0^{\tau}\omega_0(I_0, Y(\tau_1, I_0), X(\tau_1, I_0))\,\dif\tau_1$$
is a strip around the real axis that contains a domain
$$D_\xi=\bigg\{\xi\colon \Big\{ \begin{array} {lc}
|\mathrm{Im}\;\xi|<\sigma_1+\delta_1|\mathrm{Re}\;\xi|, & |\mathrm{Re}\;\xi|>c_5\Gamma \\
|\mathrm{Im}\;\xi|<\sigma_2, & |\mathrm{Re}\;\xi|\leq c_5\Gamma
\end{array}  \bigg\}$$
\begin{figure}[H]
\centering
\includegraphics[width=10cm]{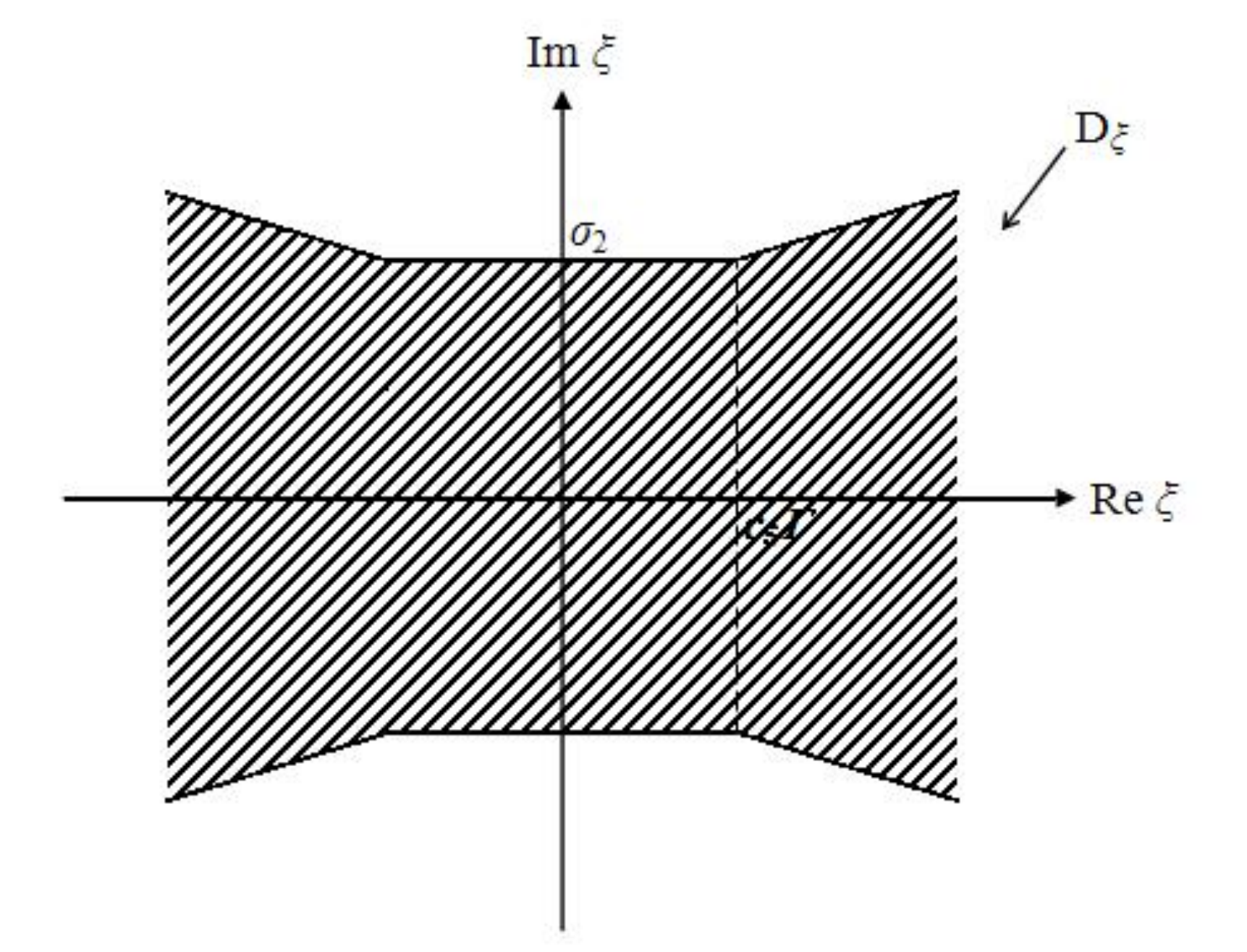}
\caption{}
\end{figure}

\noindent where $\sigma_1$, $\sigma_2$, $\delta_1$, $c_5$ are positive constants and $\sigma_2=\sigma_1+\delta_1\cdot c_5\Gamma$ for continuity. Straight lines $\mathrm{Im}\,\xi=\pm\gamma$ lie in $D_\xi$ and have a positive distance from the boundary of $D_\xi$ (i.e. $\gamma<\sigma_2$).
\end{lem}

\begin{lem}
The function $F_1(y,x,\ph)$ satisfies the estimate:
$$|F_1(\widehat Y(\xi,I_0), \widehat X(\xi,I_0), \ph)|<\frac{\mathrm{const}}{1+|\xi|^{2+\nu}}.$$
\end{lem}

\begin{proof}
First of all let us list two forms of Hamiltonian: \vskip 5pt
\begin{eqnarray*}
&&H_0(I, y, x)+\eps H_1(I, \ph, y,x)=h_0, \\
&&-I=F_0(y,x)+\eps F_1(y,x,\ph).
\end{eqnarray*}
Therefore, we have
\begin{eqnarray*}
&&H_0(-(F_0+\eps F_1), y, x)+\eps H_1(-(F_0+\eps F_1), \ph, y,x)=h_0. \\
\textup{Then let}\ \eps=0: &&H_0(-F_0,y,x)=h_0.
\end{eqnarray*}
So
\begin{eqnarray*}
&&H_0(-(F_0+\eps F_1), y, x)+\eps H_1(-(F_0+\eps F_1), \ph, y,x)=H_0(-F_0,y,x), \\
&&H_0(-(F_0+\eps F_1), y, x)-H_0(-F_0,y,x)=-\eps H_1(-(F_0+\eps F_1), \ph, y,x).
\end{eqnarray*}
Applying the intermediate value theorem here, we get that there exist some value $\widetilde I$, such that
\begin{eqnarray*}
&&\frac{\partial H_0(\widetilde I,y,x)}{\partial I}\cdot(-\eps F_1)=-\eps H_1(-(F_0+\eps F_1), \ph, y,x) \\
\textup{so}&&F_1(y,x,\ph)=\frac1{\displaystyle\frac{\partial H_0(\widetilde I,y,x)}{\partial I}}\cdot H_1(I, \ph, y,x)
\end{eqnarray*}
We know from assumptions that
\begin{eqnarray*}
&&|\omega_0|=\left|\frac{\partial H_0}{\partial I}\right|>\mathrm{const} \\
&\textup{and}&|H_1(I, \ph, Y(\tau, I_0), X(\tau, I_0))|<\frac c{1+\left|\int\limits_0^{\tau}\omega_0(I_0, Y(\tau_1, I_0), X(\tau_1, I_0))\,\dif\tau_1\right|^{2+\nu}}
\end{eqnarray*}
Therefore,
\begin{eqnarray*}
|F_1(\widehat Y(\xi,I_0), \widehat X(\xi,I_0), \ph)|&<&\mathrm{const}\cdot|H_1(I, \ph, Y(\tau, I_0), X(\tau, I_0))| \\
&<&\frac{\mathrm{const}}{1+|\xi|^{2+\nu}}.
\end{eqnarray*}

\end{proof}

\begin{cor}
For $\ph\in\widetilde D_\ph=D_\ph-\delta_\ph$, $\xi\in\widetilde D_\xi=D_\xi-\delta_\xi$,
$$\left|\frac{\partial F_1(\widehat Y(\xi,I_0), \widehat X(\xi,I_0), \ph, h_0)}{\partial\ph}\right|<\frac{\mathrm{const}}{1+|\xi|^{2+\nu}},$$
$$\left|\frac{\partial F_1(\widehat Y(\xi,I_0), \widehat X(\xi,I_0), \ph, h_0)}{\partial\xi}\right|<\frac{\mathrm{const}}{1+|\xi|^{2+\nu}},$$
where $\delta_\xi$ is a positive constant.
\end{cor}

Now let us draw the phase portrait of $F_0(y,x)=-I$:
\begin{figure}[H]
\centering
\includegraphics[width=11cm]{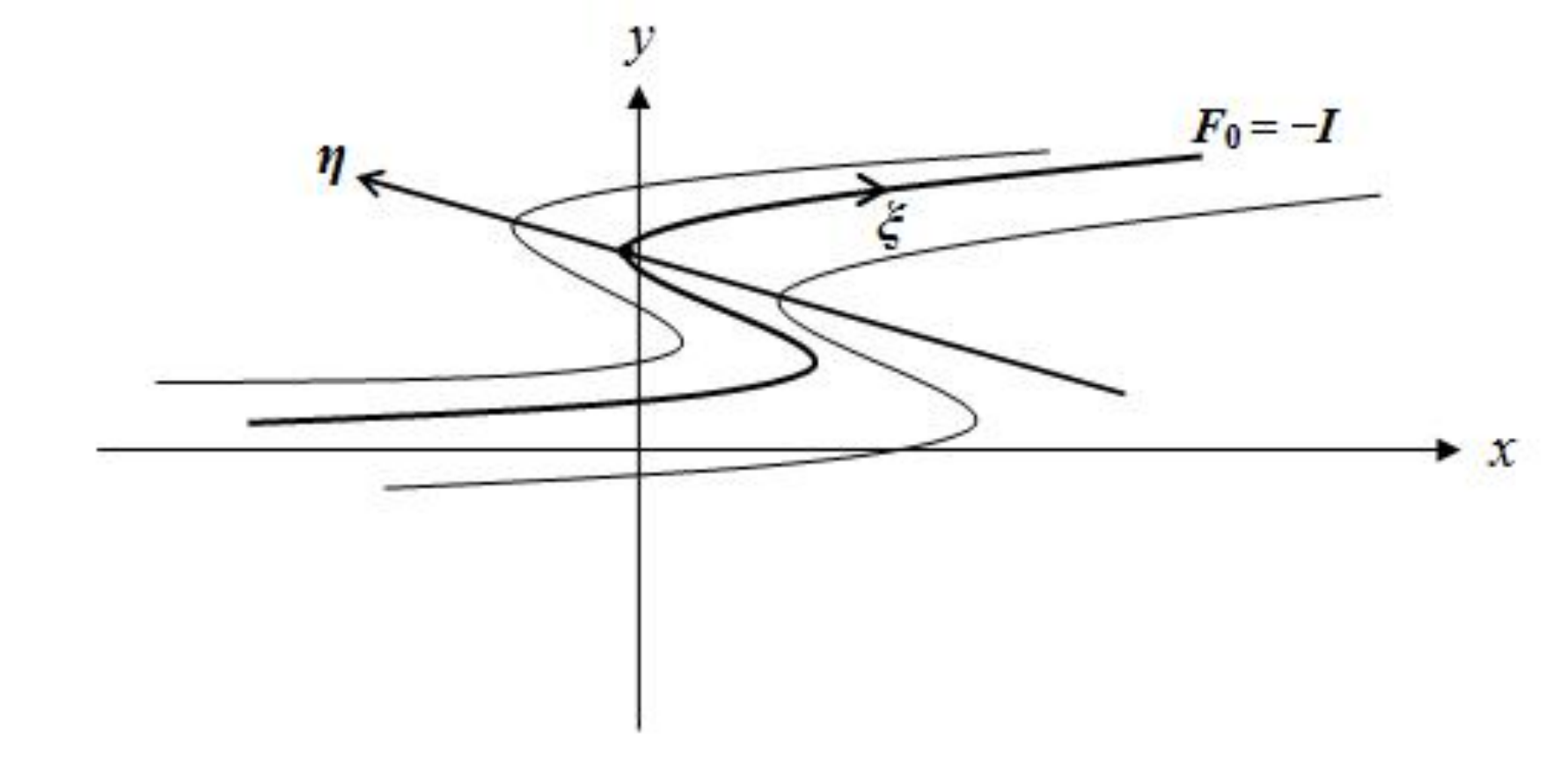}
\caption{}
\end{figure}

\noindent and introduce new coordinates $(\xi,\eta)$. Here $\eta$ is exactly the value of $F_0$ on level lines $F_0=\eta$, and $\xi$ is the slow time of motion. The initial section is $\widehat X(0,I)$, $\widehat Y(0,I)$ parametrized by $I$. Then construct a transformation $(\xi,\eta)\mapsto(x,y)$.

\begin{lem}
In the domain $D_\xi\times D_\eta $, where $D_\xi$ is what mentioned above and  $D_\eta$ is a neighbourhood of the real point $-I_*$, the mapping $(\xi,\eta)\mapsto(x,y)$ is analytic and canonical.
\end{lem}

\begin{proof}
With map $(\xi,\eta)\mapsto(x,y)$, consider $x$ and $y$ as functions of $\xi$, $\eta$:
$$x=x(\xi,\eta),\quad y=y(\xi,\eta).$$
From the motion above, we have obtained the relations:
$$\frac{\partial x}{\partial \xi}=\frac{\partial F_0}{\partial y},\quad \frac{\partial y}{\partial \xi}=-\frac{\partial F_0}{\partial x}.$$
Let us differentiate both sides of the Hamiltonian $F_0(y,x,h_0)=\eta$ with respect to $\eta$:
$$\frac{\partial F_0}{\partial y}\frac{\partial y}{\partial\eta}+\frac{\partial F_0}{\partial x}\frac{\partial x}{\partial\eta}=1.$$
Replacing with corresponding derivatives, we get
$$\frac{\partial x}{\partial\xi}\frac{\partial y}{\partial\eta}-\frac{\partial y}{\partial\xi}\frac{\partial x}{\partial\eta}=1.$$
That is indeed the Jacobi determinant:
$$\frac{D(x,y)}{D(\xi,\eta)}=\det
\left|\begin{array}{cc}
\dfrac{\partial x}{\partial\xi} & \dfrac{\partial x}{\partial\eta} \\ \\
\dfrac{\partial y}{\partial\xi} & \dfrac{\partial y}{\partial\eta}
\end{array}\right|
=1.$$
Therefore, $(x,y)\mapsto(\xi,\eta)$ is canonical transformation.

\end{proof}

After this transformation, the new Hamiltonian has the form
$$-I(\eta,\xi,\ph)=\eta+\eps G_1(\eta,\xi,\ph).$$
Function $G_1$ is analytic and has the following estimate:
$$|G_1(\eta,\xi,\ph)|<\frac{\mathrm{const}}{1+|\xi|^{2+\nu}}.$$
Also from Cauchy estimate \cite{8}, for $\eta\in\widetilde D_\eta=D_\eta-\delta_\eta$, $\xi\in\widetilde D_\xi=D_\xi-\delta_\xi$, $\ph\in\widetilde D_\ph$, where $\delta_\eta$ and $\delta_\xi$ are positive constants, we obtain
$$\left|\frac{\partial G_1(\eta,\xi,\ph)}{\partial\ph}\right|<\frac{\mathrm{const}}{1+|\xi|^{2+\nu}},$$
$$\left|\frac{\partial G_1(\eta,\xi,\ph)}{\partial\eta}\right|<\frac{\mathrm{const}}{1+|\xi|^{2+\nu}},$$
$$\left|\frac{\partial G_1(\eta,\xi,\ph)}{\partial\xi}\right|<\frac{\mathrm{const}}{1+|\xi|^{2+\nu}}.$$

Because of the above estimate, we can take another similar isoenergetic reduction through expressing $\eta$, via $I$, $\ph$, $\xi$ and constant $h_0$:
$$-\eta(I,\ph,\xi)=I+\eps K_1(I,\ph,\xi).$$
Here $(-\eta)$ is the new Hamiltonian, and $\xi$ is the new slow time. $I$ and $\ph$ are conjugate variables.

From implicit function theorem, we can obtain the similar conclusions that function $K_1$ is analytic in $\widehat D_I\times D_\ph\times D_\xi$ and estimates of $K_1$ and its derivatives have the following forms:
\begin{eqnarray*}
&&|K_1(I,\ph,\xi)|<\frac{\mathrm{const}}{1+|\xi|^{2+\nu}}, \\
&&\left|\frac{\partial K_1(I,\ph,\xi)}{\partial I}\right|<\frac{\mathrm{const}}{1+|\xi|^{2+\nu}}, \\
&&\left|\frac{\partial K_1(I,\ph,\xi)}{\partial\ph}\right|<\frac{\mathrm{const}}{1+|\xi|^{2+\nu}}, \\
&&\left|\frac{\partial K_1(I,\ph,\xi)}{\partial\xi}\right|<\frac{\mathrm{const}}{1+|\xi|^{2+\nu}},
\end{eqnarray*}
in $(\widehat D_I-\delta_I)\times\widetilde D_\ph\times\widetilde D_\xi$. Here $\widehat D_I$ is the same as the domain $-D_\eta$.

In this latest Hamiltonian, $K_0\equiv I$ and the frequency is $\displaystyle\frac{\partial K_0}{\partial I}=1$. The level lines
$$\mathrm{Im}\int\limits_0^\xi\dif\xi_1=\mathrm{Im}\ \xi=B=\mathrm{const}, \quad 0\leq|B|\leq\gamma$$
also lie in the domain $D_\xi$ and have a positive distance from the boundary of $D_\xi$. The constant $\gamma$ is the same as introduced in Assumption $3^\circ$.

Denoting the new time as $\thet=\eps^{-1}\xi$, we can write differential equations of motion:
$$\frac{\dif I}{\dif\thet}=-\eps\frac{\partial K_1}{\partial\ph},\quad \frac{\dif\ph}{\dif\thet}=1+\frac{\partial K_1}{\partial I}.$$

Now consider the exact solution $I(\thet)$, $\ph(\thet)$ with real initial conditions $I(0)$, $\ph(0)$ at $\thet_0=\eps^{-1}\xi_0=0$. We obtain the relation which the solution satisfies as follows:
$$I(\thet)=I(0)-\eps\int\limits_0^\thet\frac{\partial K_1(I(\thet_1),\ph(\thet_1),\eps\thet_1)}{\partial\ph}\,\dif\thet_1.$$
From the Cauchy criterion \cite{4}, we can easily prove the existence of limiting values of $I(\thet)$ as $\thet\to\pm\infty$. Moreover, we also know differential equations of motion:
\begin{eqnarray*}
&&\frac{\dif\ph}{\dif\thet}=1+\eps\frac{\partial K_1}{\partial\ph} \\
\textup{and}&&\frac{\dif\ph}{\dif t}=\frac{\partial H_0}{\partial I}+\eps\frac{\partial H_1}{\partial I}.
\end{eqnarray*}
It is evident that $\ph\to\pm\infty$ as $\thet\to\pm\infty$. Meanwhile, $t\to\pm\infty$. Therefore, the limits when $t\to\pm\infty$ are equal to those when $\thet\to\pm\infty$. Let us define values
$$I_\pm=\lim\limits_{t\to\pm\infty}I(t),\quad \Delta I=I_+-I_-.$$

Now the system under consideration is described by the Hamiltonian
$$-\eta(I,\ph,\xi)=I+\eps K_1(I,\ph,\xi),$$
where $I$ and $\ph$ are conjugate canonical variables. The action-angle variables $I$, $\ph$ of the unperturbed problem are defined in region $\widehat D_I\times D_\ph$. $\xi$ is the slow time variable defined in region $D_\xi$ and time $\thet=\eps^{-1}\xi$. The function $K_1$ is an analytic, $2\pi$--periodic in $\ph$ function. The frequency of unperturbed motion is 1. We have proved the following properties:

$1^\circ$. The functions $K_0=I$, $K_1$ can be analytically extended into a complex domain $\widehat D_I\times D_\ph\times D_\xi$, where $\widehat D_I$ is some neighbourhood of a given real point $I_*$, $D_\ph$ is a strip of some fixed width about the real axis, and $D_\xi$ is a strip
$$D_\xi=\bigg\{\xi\colon \Big\{ \begin{array} {lc}
|\mathrm{Im}\;\xi|<\sigma_1+\delta_1|\mathrm{Re}\;\xi|, & |\mathrm{Re}\;\xi|>c_5\Gamma \\
|\mathrm{Im}\;\xi|<\sigma_2, & |\mathrm{Re}\;\xi|\leq c_5\Gamma
\end{array}  \bigg\}$$ where $\sigma_1+\delta_1\cdot c_5\Gamma=\sigma_2$.
The function $K_1$ satisfies the estimate $$|K_1(I,\ph,\xi)|<\frac C{1+|\xi|^{2+\nu}}.$$ Here $\sigma$, $\delta$, $C$, $\nu$ are positive constant.

$2^\circ$. The level lines
$$\mathrm{Im}\int\limits_0^\xi\dif\xi_1=\mathrm{Im}\ \xi=B=\mathrm{const}, \quad 0\leq|B|\leq\gamma$$
lie in the domain $D_\xi$ and have a positive distance from the boundary of $D_\xi$.

Consider a solution $I(\thet)$, $\ph(\thet)$ of the Hamiltonian system with real initial conditions $I(0)$, $\ph(0)$ at $\thet=\thet_0=\eps^{-1}\xi_0=0$. For $$I_\pm=\lim\limits_{\thet\to\pm\infty}I(\thet),\quad \Delta I=I_+-I_-,$$
referring to the paper \cite{1}, we have the estimate of $\Delta I$ that is exponentially small:
$$\Delta I=O(\me^{-\frac\gamma\eps}).$$

Actually, all variables $I$, $\ph$, $y$ and $x$ mentioned before should contain the bar symbol
``$\,\bar\ \bar\ \,$". We obtain in fact the conclusion that
$$\Delta\bar I=\bar I_+-\bar I_-=O(\me^{-\frac\gamma\eps}),\quad\gamma=\mathrm{const}>0.$$

Now return the bar. Recall that from the original Hamiltonian $E(p,q,y,x)$, we transform to two forms of Hamiltonian. With generating function $S(I,q,y,x)$ where $y$, $x$ are considered as parameters, we obtain $H_0(I,y,x)$ and with $\eps^{-1}\bar yx+S(\bar I,q,\bar y,x)$ we obtain $H_0(\bar I,\bar y,\bar x)+\eps H_1(\bar I,\bar\ph,\bar y,\bar x)$. Since the limiting values of $\bar I$ exist, now let us prove the lemma that the limiting value of $I$ is equal to that of $\bar I$ and thus exists.

\begin{lem}
$$\lim\limits_{t\to\pm\infty}I(t)=\lim\limits_{t\to\pm\infty}\bar I(t).$$
\end{lem}

\begin{proof}
With generating function $S(p,q,y,x)$, the Hamiltonian $E(p,q,y,x)$ is canonically transformed to $H_0(I,y,x)$ and $p=\dfrac{\partial S(I,q,y,x)}{\partial q}$.

As we know $$\frac{\partial S}{\partial x}=\int\limits_0^t\left(\left<\frac{\partial E}{\partial x}\right>-\frac{\partial E}{\partial x}\right)\,\dif t_1,$$
so $$\frac{\partial S}{\partial x}\to0,\quad \mathrm{since}\ \frac{\partial E}{\partial x}\to0,\ \mathrm{as}\ x\to\pm\infty.$$

On the other hand, with generating function $\eps^{-1}\bar yx+S(\bar I,q,\bar y,x)$, the Hamiltonian $E(p,q,y,x)$ is transformed to $H_0(\bar I,\bar y,\bar x)+\eps H_1(\bar I,\bar\ph,\bar y,\bar x)$ and $p=\dfrac{\partial S(\bar I,q,\bar y,x)}{\partial q}$. Also we know
$$H_1=-\left.\frac{\partial H_0}{\partial\bar x}\right|_{\widetilde x}\frac{\partial S}{\partial\bar y}+\left.\frac{\partial E}{\partial y}\right|_{\widetilde y}\frac{\partial S}{\partial x},$$
where $\widetilde x$ and $\widetilde y$ are some intermediate values. As $\dfrac{\partial H_0}{\partial \bar x}\to0$, $\dfrac{\partial S}{\partial x}\to0$, thus $H_1\to0$, as $x\to\pm\infty$.

From assumption $\Big(\dfrac{\partial H_0}{\partial\bar y}\Big)^2+\Big(\dfrac{\partial H_0}{\partial\bar x}\Big)^2>\mathrm{const}$, and $\dfrac{\partial H_0}{\partial\bar x}\to0$, as $\bar x\to\pm\infty$, we know that $\left|\dfrac{\partial H_0}{\partial\bar y}\right|>\const$ for big $|\bar x|$. Thus, when $|\bar x|>c_6$, $\left|\dfrac{\partial H_0}{\partial\bar y}\right|>c_7^{-1}$, where $c_6$, $c_7$ are positive constants. From $\dot{\bar y}=\dfrac{\dif\bar y}{\dif t}=-\eps\left(\dfrac{\partial H_0}{\partial\bar x}+\eps\dfrac{\partial H_1}{\partial\bar x}\right)$, $\dot{\bar x}=\dfrac{\dif\bar x}{\dif t}=\eps\left(\dfrac{\partial H_0}{\partial\bar y}+\eps\dfrac{\partial H_1}{\partial\bar y}\right)$, the phase portrait of $H_0+\eps H_1=\const$ in $(\bar x,\bar y)$ plane will close to that of the adiabatic approximation $H_0=\const$ for $|\bar x|\leq c_6+1$.

\begin{figure}[H]
\centering
\includegraphics[width=9cm]{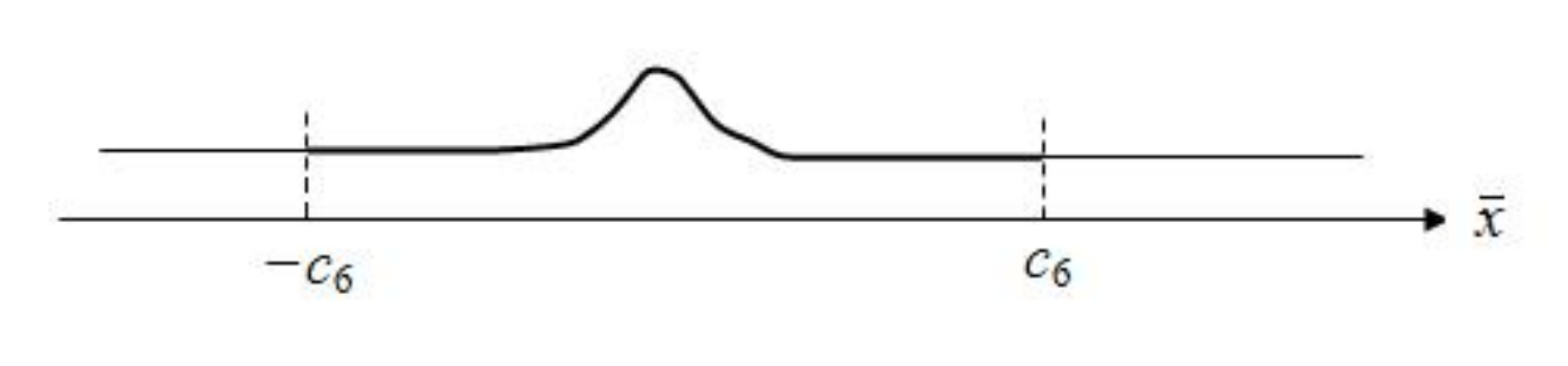}
\caption{}
\end{figure}

Thus for some time moments $t_+$ and $t_-$, we have
$$\bar x(t_+)>c_6,\quad \bar x(t_-)<-c_6.$$
Then for $t>t_+$, we have
$$\dot{\bar x}=\eps\left(\frac{\partial H_0}{\partial\bar y}+\eps\frac{\partial H_1}{\partial\bar y}\right)>\eps\big(c_7^{-1}+O(\eps)\big)>\frac12\,\eps\,c_7^{-1}.$$
Therefore, $\bar x(t)\to+\infty$ as $t\to+\infty$. Similarly, $\bar x(t)\to-\infty$ as $t\to-\infty$.

From $$y=\bar y+\eps \frac{\partial S(\bar I, q, \bar y, x)}{\partial x}=\bar y+\eps V(\bar I,\bar\ph,\bar y,x),$$
taking limiting values on both sides as $t\to\pm\infty$, we obtain $y-\bar y\to0$.

We also know that
\begin{eqnarray*}
p=P(I,\ph,y,x)=P(\bar I,\bar\ph,\bar y,x), \\
q=Q(I,\ph,y,x)=Q(\bar I,\bar\ph,\bar y,x).
\end{eqnarray*}
As  $y-\bar y\to0$, it is evident that $I-\bar I\to0$, $\ph-\bar\ph\to0$ as $t\to\pm\infty$.

Therefore, $$\lim\limits_{t\to\pm\infty}I(t)=\lim\limits_{t\to\pm\infty}\bar I(t).$$

\end{proof}

Thus $I$ and $\bar I$ have the same limit, and $\Delta I=\Delta\bar I$. Therefore, we have the following corollary:

\begin{cor}
The estimate $$\Delta I=O(\me^{-\frac\gamma\eps})$$ is valid.
\end{cor}

\bigskip
\section*{Acknowledgements}
I am glad to express my gratitude to Professor Anatoly Neishtadt for providing detailed comments and frequent discussion during the research and completion of this article.


\bigskip
\begin{thebibliography} {[99]}
\footnotesize

\bibitem{2} Arnold, V.I., Kozlov, V.V. and Neishtadt, A.I.: \emph{Mathematical Aspects of Classical and Celestial Mechanics}. Second augmented and revised edition. Springer-Verlag, 2006.
\bibitem{1} Neishtadt, A.I.: \emph{On the Accuracy of Persistence of Adiabatic Invariant in Single-Frequency Systems}. Regular and Chaotic Dynamics, Vol. 5, No 2, 2000, pp. 213-218.
\bibitem{3} Arnold, V.I.: \emph{Mathematical Methods of Classical Mechanics}. Second Edition. Springer, 2000.
\bibitem{9} Neishtadt, A.I.: \emph{On the Change in the Adiabatic Invariant on Crossing a Separatrix in System with Two Degrees of Freedom}. PMM U.S.S.R, Vol. 51, No. 5, 1987, pp. 586-592.
\bibitem{8} Brown, J.W. and Churchill, R.V.: \emph{Complex Variables and Applications}. Sixth Edition. McGraw-Hill, 1996.
\bibitem{4} Arnold, V.I.: \emph{Ordinary Differential Equations}. The MIT Press, Cambridge, Massachusetts and London, 1973.
\bibitem{5} Benettin, G., Carati, A. and Gallavotti, G.: \emph{A Rigorous Implementation of the Jeans-Landau-Teller Approximation for adiabatic invariants}. Nonlinearity 10, 1997, pp. 479-505.
\bibitem{10} Treschev, D.V.: \emph{The Method of Continuous Averaging in the Problem of Fast and Slow Motions}. Regular and Chaotic Dynamics, Vol. 12, No. 3/4, 1997, pp. 9-20


\end{thebibliography}
\end{document}